\documentclass[11pt]{amsart}
\usepackage[T1]{fontenc}
\usepackage[utf8]{inputenc}
\usepackage{lmodern}
\usepackage{amsmath,amssymb,amsthm,mathtools}
\usepackage{microtype}
\usepackage{graphicx}
\usepackage{tikz}
\usetikzlibrary{arrows.meta}
\usepackage{needspace}
\usepackage{cite}
\usepackage[hidelinks]{hyperref}
\numberwithin{equation}{section}
\newtheorem{theorem}{Theorem}[section]
\newtheorem{corollary}[theorem]{Corollary}

\theoremstyle{definition}
\newtheorem{example}[theorem]{Example}
\theoremstyle{remark}
\newtheorem{remark}[theorem]{Remark}
\DeclareMathOperator{\sgn}{sgn}
\DeclareMathOperator{\tr}{tr}
\DeclareMathOperator{\adj}{adj}
\newcommand{\T}{\mathbb S^1}
\newcommand{\D}{\mathbb D}
\newcommand{\Lcal}{\mathcal L}
\newcommand{\Pcal}{\mathcal P}
\newcommand{\dmu}{\,\mathrm d\mu}
\newcommand{\dm}{\,\mathrm dm}
\newcommand{\dd}{\,\mathrm d}
\newcommand{\qbinom}[2]{\genfrac{[}{]}{0pt}{}{#1}{#2}_q}
\allowdisplaybreaks[1]

\title[Extensions of certain Markov inequalities]{Extensions of certain Markov inequalities to Toeplitz matrices}
\author{K. Castillo}
\address{CMUC, Department of Mathematics, University of Coimbra, 3000-143 Coimbra, Portugal}
\email{kenier@mat.uc.pt}
\author{A. Suzuki}
\address{CMUC, Department of Mathematics, University of Coimbra, 3000-143 Coimbra, Portugal}
\email{asuzuki@uc.pt}
\date{}

\begin{document}
\begin{abstract}
We establish finite Toeplitz analogues of Markov's inequalities for Hankel determinants and zeros of orthogonal polynomials. We characterise the real symmetric positive-definite Toeplitz matrices whose determinants dominate those of all such matrices below them in the alternating moment order: non-negativity of the diagonal sums of deleted minors is necessary and sufficient. Equality holds only when the matrices coincide, and a quantitative estimate controls the determinant gap. Paraorthogonal zeros in the right half-plane provide a geometric sufficient condition in both parities. For differentiable finite Hermitian moment data, we derive individual and relative angular-velocity identities. For real symmetric data, these identities yield a finite comparison of all upper zeros under simultaneous changes in the moments: a semicircle condition at one endpoint is inherited by the other, and every angular inequality is strict when the matrices differ. No regularity assumption on a representing measure is required. Coefficient criteria also establish nearest-zero monotonicity. An exact Rogers--Szeg\H{o} counterexample shows that this need not extend to every upper zero.
\end{abstract}
\maketitle

\section{Introduction}

Let $(a_j)_{j\geq0}$ be a sequence of real numbers, and write
\[
 H_m=(a_{r+s})_{r,s=0}^{m-1},\quad
 H_m^{[1]}=(a_{r+s+1})_{r,s=0}^{m-1}.
\]
The Stieltjes moment problem asks for a positive Borel measure $\mu$ on $[0,\infty)$ such that
\[
 a_j=\int_0^\infty x^j\dmu(x),\quad j\geq0.
\]
Such a representation exists if and only if $H_m$ and $H_m^{[1]}$ are positive semidefinite for every $m\geq1$. Their positive definiteness for every $m$ is the corresponding infinite-support condition. In this latter case the monic orthogonal polynomials $p_m$ exist in every degree and have simple positive zeros. The associated infinite Hankel matrix is strictly totally positive: every finite minor is positive. These classical connections between moment problems, continued fractions and matrix positivity are discussed in \cite{M94,S94A,ST43,KN77,P10}. Throughout this paper, \emph{totally non-negative} means that every minor is non-negative, whereas \emph{strictly totally positive} means that every minor is positive.

Markov considered the effect of varying a single moment while preserving the relevant Stieltjes positivity conditions. For $p_{n+1}$ and $0\leq k\leq2n$, the condition
\[
 (-1)^k a_k'(t)<0
\]
implies that every zero of $p_{n+1}(\cdot;t)$ increases and that $\det H_{n+1}(t)$ decreases; see \cite[pp.~15--21]{M94} and \cite[p.~114]{ST43}. The range matters: $p_{n+1}$ also depends on $a_{2n+1}$, whereas $\det H_{n+1}$ does not. Markov's paper appeared in 1894, and Shohat's English translation followed in 1940. These results also appear in \cite{KN77,G59}.

Our aim is to address the corresponding questions for Toeplitz matrices. Orthogonal polynomials on the unit circle and their Christoffel--Darboux kernels take the place of the real-line objects, whilst paraorthogonal polynomials supply the zeros on the circle whose motion we study. Related results on paraorthogonal zeros and their variation appear in \cite{S07,CP18,C19,C23}; general background may be found in \cite{S05I,S11}. Two-sided P\'olya frequency sequences provide a natural class of examples, but the Toeplitz positivity in their definition must be distinguished from the Hankel positivity of a Stieltjes moment sequence.

The determinant comparison below involves only finite matrices and permits simultaneous changes in several moments. Its hypothesis is necessary and sufficient for comparison throughout the alternating moment order, and its equality case yields strict monotonicity without strict total positivity. For zero motion, the essential object is a finite Hermitian moment functional: differentiation on Laurent polynomials gives an exact two-zero identity. Neither a differentiable density nor an absolutely continuous representing measure is needed. The resulting sign criteria lead to a finite comparison theorem with a semicircle hypothesis at only one endpoint; the other endpoint inherits the same geometry. Coefficient positivity also controls suitable extreme zeros without a semicircle restriction, but the final example shows that this conclusion does not extend to every zero in general.

Section~\ref{sec:results} introduces the finite moment data and presents the results. The corollaries and the one-endpoint comparison are proved immediately after their statements, and the Rogers--Szeg\H{o} example is worked out in full. The proofs of Theorems~\ref{thm:comparison}, \ref{thm:hurwitz} and~\ref{thm:variation}, together with the kernel coefficient identities, are given in Section~\ref{sec:proofs}.

\section{Finite moment data and main results}\label{sec:results}

Write $\T=\{z\in\mathbb C:|z|=1\}$, $\D=\{z\in\mathbb C:|z|<1\}$ and
\[
 \T_\pm=\{z\in\T:\ \pm\Re z>0\}.
\]
For a matrix $M$ with row and column indices starting at $0$, let $M^{(r,s)}$ denote the submatrix obtained by deleting row $r$ and column $s$. Its determinant is an \emph{unsigned deleted minor}; the corresponding cofactor is $(-1)^{r+s}\det M^{(r,s)}$. A real sequence $(c_j)_{j\in\mathbb Z}$ is a two-sided P\'olya frequency sequence if its infinite Toeplitz matrix $(c_{r-s})_{r,s\in\mathbb Z}$ is totally non-negative, in the classical convention of \cite{E52,K68}.

Fix $n\geq1$ and finite Hermitian data $c_{-j}=\overline{c_j}$, $0\leq j\leq n$. We use the convention
\begin{equation}\label{eq:gram}
 C_{m+1}=(c_{r-s})_{r,s=0}^{m},\quad 0\leq m\leq n,
 \quad C_{n+1}>0.
\end{equation}
Define the Hermitian functional on $\Lambda_n=\operatorname{span}\{z^j:-n\leq j\leq n\}$ by
\begin{equation}\label{eq:functional}
 \Lcal[z^j]=c_{-j},\quad
 \langle f,g\rangle=\Lcal[\overline f g],\quad f,g\in\Pcal_n,
\end{equation}
where $\Pcal_n$ is the space of polynomials of degree at most $n$ and conjugation is taken on $\T$. Thus $C_{n+1}$ is the Gram matrix of $1,z,\ldots,z^n$. A positive representing measure $\sigma$ exists, with
\[
 c_j=\int_{\T} z^{-j}\dd\sigma(z),\quad |j|\leq n;
\]
see \cite[Section~1.3]{S05I}. None of the results below depends on the choice or regularity of a representing measure.

Let $Q_m$ be the monic polynomial of degree $m$ orthogonal to $1,z,\ldots,z^{m-1}$, and let $q_m=\kappa_m Q_m$, where $\kappa_m=\Lcal[|Q_m|^2]^{-1/2}>0$. In particular, $Q_0=1$, and
\begin{equation}\label{eq:opuc}
 Q_m(z)=\frac{1}{\det C_m}
 \det\begin{pmatrix}
 c_0&c_{-1}&\cdots&c_{-m}\\[7pt]
 c_1&c_0&\cdots&c_{1-m}\\[7pt]
 \vdots&\vdots&&\vdots\\[7pt]
 c_{m-1}&c_{m-2}&\cdots&c_{-1}\\[7pt]
 1&z&\cdots&z^m
 \end{pmatrix},\quad 1\leq m\leq n.
\end{equation}
For $w,z\in\mathbb C$ and $b\in\T$, put
\begin{align}
 K_n(w,z)&=\sum_{m=0}^n\overline{q_m(w)}q_m(z),\label{eq:kernel}\\[7pt]
 P_{n+1}(z,b)&=zQ_n(z)-\overline b\,Q_n^*(z),
 \quad Q_n^*(z)=z^n\overline{Q_n(1/\overline z)}.\label{eq:popuc}
\end{align}
Here $K_n$ is the Christoffel--Darboux (CD) kernel, and $Q_n^*$ is the reversed polynomial, whose value at $z=0$ is defined by polynomial continuation. The polynomial $P_{n+1}$ is monic and has $n+1$ simple zeros on $\T$; the zeros of $Q_n$ lie in $\D$. We shall use these standard finite-dimensional consequences of orthogonality, together with the CD formula; see \cite{S05I,S11}. A polynomial is called \emph{Hurwitz} if all its zeros lie in the open left half-plane.

\Needspace{9\baselineskip}
\subsection{Determinant comparison}

\begin{theorem}[An exact finite comparison criterion]\label{thm:comparison}
Let $n\geq1$, and let $B=(b_{r-s})_{r,s=0}^n$ be a real symmetric positive-definite Toeplitz matrix, with $b_{-j}=b_j$. Define
\begin{equation}\label{eq:diagonalsums}
 S_k(B)=\sum_{\substack{0\leq r,s\leq n\\[7pt]|r-s|=k}}
 \det B^{(r,s)},\quad 0\leq k\leq n,
\end{equation}
where $B^{(r,s)}$ is obtained by deleting row $r$ and column $s$. The following statements are equivalent:
\begin{enumerate}
\item[(i)] $S_k(B)\geq0$ for every $0\leq k\leq n$.
\item[(ii)] Every real symmetric positive-definite Toeplitz matrix
$A=(a_{r-s})_{r,s=0}^n$, with $a_{-j}=a_j$, satisfying
\begin{equation}\label{eq:momentorder}
 (-1)^j(b_j-a_j)\geq0,\quad 0\leq j\leq n,
\end{equation}
satisfies $\det A\leq\det B$.
\end{enumerate}
Whenever these conditions hold, equality in~\textup{(ii)} occurs if and only if $A=B$. The condition $S_0(B)>0$ is automatic. In particular,~\textup{(i)} holds if all unsigned deleted minors of $B$ are non-negative, and hence whenever $B$ is totally non-negative.
\end{theorem}

The diagonal-sum condition is strictly weaker than requiring each deleted minor to be non-negative. For example,
\[
 B=\begin{pmatrix}
 4&-1&-3&3\\[7pt]-1&4&-1&-3\\[7pt]-3&-1&4&-1\\[7pt]3&-3&-1&4
 \end{pmatrix}
\]
has leading principal determinants $4,15,14,5$, and is therefore positive definite. Its diagonal sums are
\[
 (S_0(B),S_1(B),S_2(B),S_3(B))=(40,10,28,22),
\]
whereas $\det B^{(1,2)}=-1$. Thus Theorem~\ref{thm:comparison} applies beyond the class defined by non-negativity of all deleted minors.

Only $B$ is required to satisfy the diagonal-sum condition, and no differentiability is assumed. Consequently, if $C_{n+1}(t)=(c_{r-s}(t))_{r,s=0}^n$ is real symmetric and positive definite, $S_k(C_{n+1}(t))\geq0$ for $0\leq k\leq n$, and each $(-1)^jc_j(t)$ is non-decreasing, then its determinant is non-decreasing. It is strictly increasing if, between any two parameter values, at least one signed moment increases strictly. For a single varying pair $c_{\pm k}$, $0\leq k\leq n$, this yields the precise comparison
\begin{equation}\label{eq:singlecomparison}
 \sgn\!\left(\det C_{n+1}(t_2)-\det C_{n+1}(t_1)\right)
 =\sgn\!\left((-1)^k(c_k(t_2)-c_k(t_1))\right).
\end{equation}
Here the family need not be monotone: the theorem is applied with the larger signed moment as the second endpoint. Finite total non-negativity suffices in all these statements; an infinite P\'olya frequency extension is unnecessary.

The next theorem gives a geometric sufficient condition for the diagonal sums by establishing positivity of every deleted minor. The boundary parameter depends on the parity of $n$.

\Needspace{9\baselineskip}
\begin{theorem}[A half-plane criterion]\label{thm:hurwitz}
Let $n\geq1$ and $c_{-j}=c_j\in\mathbb R$ for $0\leq j\leq n$, and suppose that $C_{n+1}=(c_{r-s})_{r,s=0}^n>0$. Let $Q_n$ be its monic orthogonal polynomial, defined by \eqref{eq:opuc}, and put
\[
 P(z)=zQ_n(z)-(-1)^nQ_n^*(z).
\]
If every zero of $P$ belongs to $\T_+$, or equivalently $P(-z)$ is Hurwitz, then
\[
 \det C_{n+1}^{(r,s)}>0,\quad 0\leq r,s\leq n.
\]
Thus $C_{n+1}$ may play the role of $B$ in Theorem~\ref{thm:comparison}.
\end{theorem}

When $n$ is even, $P=P_{n+1}(\cdot,1)$ has the fixed zero $1$ and is a non-zero constant multiple of $(z-1)K_n(1,z)$. Hence the hypothesis in Theorem~\ref{thm:hurwitz} is precisely that $K_n(1,-z)$ be Hurwitz. For odd $n$, the appropriate polynomial is $P_{n+1}(\cdot,-1)$; a Hurwitz condition on $K_n(1,-z)$ would be impossible, because $P_{n+1}(\cdot,1)$ also has the zero $-1$. The half-plane formulation therefore covers both parities.

\subsection{Exact identities for angular velocities}

In the following theorem all the finite moments may vary. A dot on the functional denotes differentiation of the moments, with the Laurent-polynomial argument held fixed.

\Needspace{9\baselineskip}
\begin{theorem}[Angular variation]\label{thm:variation}
Let $n\geq1$, let $J$ be an open real interval, and let $c_0,\ldots,c_n$ be $C^1$ functions on $J$, with $c_0$ real-valued and $c_{-j}=\overline{c_j}$. Suppose that $C_{n+1}(t)=(c_{r-s}(t))_{r,s=0}^n$ is positive definite for every $t\in J$. Let $Q_n(\cdot;t)$ be the monic orthogonal polynomial defined by \eqref{eq:opuc}, let $b:J\to\T$ be $C^1$, and set
\[
 P(z;t)=zQ_n(z;t)-\overline{b(t)}Q_n^*(z;t),\quad
 \Lcal_t[z^j]=c_{-j}(t),\quad |j|\leq n.
\]
Every zero of $P$ has a $C^1$ parametrisation on $J$ and a real $C^1$ argument. Write $b=e^{i\beta}$ with a real $C^1$ function $\beta$, and set $\nu=\Lcal_t[|Q_n|^2]=\kappa_n^{-2}$. For a zero branch $\xi=e^{i\theta_\xi}$, define
\begin{align}
 h_\xi(t)&=\Lcal_t\!\left[\left|\frac{P(z;t)}{z-\xi(t)}\right|^2\right]>0,\label{eq:h}\\[7pt]
 V_\xi(z;t)&=2\nu(t)\Im\left(
 \frac{Q_n(z;t)K_{n-1}(z,\xi(t);t)}{Q_n(\xi(t);t)}\right),
 \quad z\in\T,\label{eq:V}
\end{align}
where $K_{n-1}(w,z;t)=\sum_{m=0}^{n-1}\overline{q_m(w;t)}q_m(z;t)$ is the CD kernel formed from the orthonormal polynomials $q_m=\kappa_mQ_m$ for the same finite moment data. Then $V_\xi\in\Lambda_n$ is real-valued on $\T$, and
\begin{equation}\label{eq:individual}
 h_\xi\theta_\xi'=-\nu\beta'+\dot\Lcal_t[V_\xi].
\end{equation}
For two distinct branches $\zeta=e^{i\varphi}$ and $\eta=e^{i\psi}$, define
\begin{align}
 B_{\zeta,\eta}(z;t)&=
 i(\zeta(t)-\eta(t))\frac{z|P(z;t)|^2}
 {(z-\zeta(t))(z-\eta(t))},\quad z\in\T.\label{eq:B}
\end{align}
Then $B_{\zeta,\eta}$ is a real-valued Laurent polynomial in $\Lambda_n$, and
\begin{equation}\label{eq:twozero}
 h_\zeta\varphi'-h_\eta\psi'=\dot\Lcal_t[B_{\zeta,\eta}].
\end{equation}
For a real-valued Laurent polynomial $F=\sum_{j=-n}^n f_jz^j$ on $\T$, the derivative on the right is explicitly
\begin{equation}\label{eq:D}
 \dot\Lcal_t[F]=f_0c_0'(t)+2\Re\sum_{j=1}^n\overline{f_j}\,c_j'(t).
\end{equation}
In particular, if only $c_{\pm k}$ vary, where $0\leq k\leq n$, then
\[
 \dot\Lcal_t[F]=(2-\delta_{k0})\Re(\overline{f_k}\,c_k').
\]
Here $\delta_{k0}=1$ if $k=0$ and $0$ otherwise.
\end{theorem}

To obtain monotonicity criteria, we use either a conjugate pair of moving zeros or one moving zero and a fixed zero. In both cases, the direction is determined by explicit Laurent coefficients. We use the convention $\sgn0=0$.

\begin{remark}[Conjugate and fixed zeros]\label{rem:signs}
Under the hypotheses of Theorem~\ref{thm:variation}, suppose first that $\eta(t)=\overline{\zeta(t)}\ne\zeta(t)$ for every $t\in J$. Write
\begin{equation}\label{eq:A}
 A_\zeta(z;t)=\frac{z|P(z;t)|^2}{(z-\zeta(t))(z-\overline{\zeta(t)})}
 =\sum_{j=-n}^n d_j(t)z^j.
\end{equation}
This Laurent polynomial is real on $\T$, and \eqref{eq:twozero} becomes
\begin{equation}\label{eq:conjugate}
 (h_\zeta+h_{\overline\zeta})\varphi'
 =-2\Im\zeta\,\dot\Lcal_t[A_\zeta].
\end{equation}
If only $c_{\pm k}$ vary, then
\begin{equation}\label{eq:conjugatesign}
 \sgn\varphi'=-\sgn(\Im\zeta)\,
 \sgn\Re(\overline{d_k}c_k').
\end{equation}
If instead $\eta\in\T$ is a fixed zero distinct from $\zeta(t)$, write $B_{\zeta,\eta}(z;t)=\sum_{j=-n}^n e_j(t)z^j$. If only $c_{\pm k}$ vary, where $0\leq k\leq n$, then
\begin{equation}\label{eq:fixedsign}
 h_\zeta\varphi'=(2-\delta_{k0})\Re(\overline{e_k}c_k'),\quad
 \sgn\varphi'=\sgn\Re(\overline{e_k}c_k').
\end{equation}
Both identities remain valid for arbitrary simultaneous moment changes by using \eqref{eq:D} in place of the single-moment expression.
\end{remark}

If the moments are fixed, \eqref{eq:individual} gives $\theta_\xi'=-\nu\beta'/h_\xi$: every zero moves in the direction opposite to that of the boundary parameter $b$. More generally, the product of the zeros is $(-1)^n\overline b$, so the sum of their angular velocities is $-\beta'$. This is consistent with simultaneous moment variation and explains why the boundary-parameter term cancels in \eqref{eq:twozero}.

The coefficients also admit CD-kernel representations, which link the finite-functional argument to Fourier coefficient formulae. Let $\dm(z)=\dd\theta/(2\pi)$ for $z=e^{i\theta}$ denote normalised arc measure, used here solely to extract Laurent coefficients. In the following formulae, any omitted parameter argument is the same fixed $t$. In the conjugate case, the coefficients $d_k$ in \eqref{eq:A} satisfy
\begin{equation}\label{eq:bkconj}
\begin{aligned}
 \kappa_n^4|Q_n(\zeta)|^2d_k&=-\int_0^{2\pi}e^{-ik\theta}
 \frac{\sin((\theta-\varphi)/2)}{\sin((\theta+\varphi)/2)}
 |K_n(\zeta,e^{i\theta})|^2\,\frac{\dd\theta}{2\pi}\\[7pt]
 &=-\frac{\zeta Q_n(\zeta)}{Q_n(\overline\zeta)}
 \int_{\T} z^{-k}K_n(\zeta,z)K_n(z,\overline\zeta)\dm(z).
\end{aligned}
\end{equation}
In the fixed-zero case, the coefficients $e_k$ satisfy
\begin{equation}\label{eq:bkfixed}
 \kappa_n^4e_k=-\frac{i\overline\eta(\zeta-\eta)}
 {Q_n(\eta)\overline{Q_n(\zeta)}}
 \int_{\T} z^{-k}K_n(\zeta,z)K_n(z,\eta)\dm(z).
\end{equation}
The apparent singularities in \eqref{eq:A}, \eqref{eq:B} and \eqref{eq:bkconj} are removable; the expressions are understood through their polynomial quotients. The denominators involving $Q_n$ do not vanish on $\T$. The factors $\kappa_n^4|Q_n(\zeta)|^2$ and $\kappa_n^4$ are positive, so these rescalings leave the sign formulae unchanged. The use of $m$ for coefficient extraction imposes no condition on a measure representing the moments.

\subsection{Consequences for zero motion}

With real moment data and $b\in\{-1,1\}$, the paraorthogonal zeros occur in conjugate pairs, apart from possible real zeros. We parametrise an upper-half-plane zero $\zeta(t)=e^{i\varphi(t)}$ by $0<\varphi(t)<\pi$; a decreasing argument corresponds to clockwise motion. Conjugation reverses each direction for the corresponding lower-half-plane zero.

\Needspace{7\baselineskip}
\begin{corollary}[The right semicircle]\label{cor:right}
Let $n\geq1$, $0\leq k\leq n$, and let $J$ be an open real interval. Suppose that $C_{n+1}(t)=(c_{r-s}(t))_{r,s=0}^n$ is real symmetric and positive definite on $J$, with $C^1$ entries and only $c_{\pm k}$ varying. Let $P_{n+1}$ be defined by \eqref{eq:opuc}--\eqref{eq:popuc}, and put $b=(-1)^n$. If all its zeros belong to $\T_+$ for every $t\in J$, then every upper-half-plane zero moves clockwise whenever $(-1)^kc_k(t)$ is decreasing. The conclusion is non-strict for a non-increasing signed moment and strict for a strictly decreasing one.
\end{corollary}

\begin{proof}
Fix a conjugate pair $\zeta,\overline\zeta$, and put $R=P_{n+1}/((z-\zeta)(z-\overline\zeta))$. Its degree is $n-1$, and its roots lie in $\T_+$. Therefore $(-1)^{n-1}R(-z)$ has strictly positive coefficients, as follows by factoring it into real linear factors and real quadratic factors. If $R(z)=\sum_{j=0}^{n-1}f_jz^j$, then
\[
 (-1)^{n-1-j}f_j>0,\quad
 |R(z)|^2=\sum_{j=-n+1}^{n-1}g_jz^j,\quad (-1)^jg_j>0.
\]
On $\T$,
\begin{equation}\label{eq:Ag}
 A_\zeta(z)=(z+z^{-1}-2\Re\zeta)|R(z)|^2,
 \quad d_k=g_{k-1}+g_{k+1}-2\Re\zeta\,g_k,
\end{equation}
where $g_j=0$ outside $[-n+1,n-1]$. Every non-zero summand has sign $(-1)^{k+1}$, and at least one is non-zero for $0\leq k\leq n$. Thus $(-1)^{k+1}d_k>0$. Formula~\eqref{eq:conjugate} gives
\[
 \varphi'=a(t)\,\frac{\dd}{\dd t}\bigl((-1)^kc_k(t)\bigr),\quad
 a(t)=\frac{2(2-\delta_{k0})\Im\zeta(t)\,|d_k(t)|}
 {h_\zeta(t)+h_{\overline\zeta}(t)}>0.
\]
The coefficient $a$ is continuous and has a positive minimum on every compact subinterval. Integration therefore proves strict monotonicity even when the derivative of the strictly monotone moment vanishes at some points.
\end{proof}

\Needspace{7\baselineskip}
\begin{corollary}[The left semicircle]\label{cor:left}
Let $n\geq1$, $0\leq k\leq n$, and let $J$ be an open real interval. Suppose that $C_{n+1}(t)=(c_{r-s}(t))_{r,s=0}^n$ is real symmetric and positive definite on $J$, with $C^1$ entries and only $c_{\pm k}$ varying. Let $P_{n+1}$ be defined by \eqref{eq:opuc}--\eqref{eq:popuc}, and put $b=-1$. If all its zeros belong to $\T_-$ for every $t\in J$, then every upper-half-plane zero moves clockwise whenever $c_k(t)$ is increasing. The conclusion is non-strict or strict according as the moment is non-decreasing or strictly increasing.
\end{corollary}

\begin{proof}
Here the polynomial $R$ in the preceding proof is real monic and Hurwitz, so its coefficients are positive. Consequently all coefficients $g_j$ in their natural range are positive. Since $\Re\zeta<0$, \eqref{eq:Ag} gives $d_k>0$ for $0\leq k\leq n$. Formula~\eqref{eq:conjugate} and the same compact-interval argument prove the claim.
\end{proof}

The third corollary requires no half-plane restriction on the zeros. Its coefficient condition follows from finite total non-negativity; strict total positivity is unnecessary.

\Needspace{9\baselineskip}
\begin{corollary}[The zero nearest to $-1$]\label{cor:pf}
Let $n\geq1$, $0\leq k\leq n$, and let $J$ be an open real interval. Suppose that $C_{n+1}(t)=(c_{r-s}(t))_{r,s=0}^n$ is real symmetric and positive definite on $J$, with $C^1$ entries and only $c_{\pm k}$ varying. Let $P_{n+1}$ be defined by \eqref{eq:opuc}--\eqref{eq:popuc}, put $b=(-1)^n$, and set
\[
 p(z;t)=(-1)^{n+1}P_{n+1}(-z,b;t).
\]
Suppose that all coefficients of $p(\cdot;t)$ are non-negative for every $t\in J$. This holds, in particular, if $C_{n+1}(t)$ is totally non-negative. Let $\zeta(t)$ be the upper-half-plane zero of $P_{n+1}(\cdot,b;t)$ nearest to $-1$. For $1\leq k\leq n$, this zero moves clockwise: non-strictly when $(-1)^kc_k(t)$ is non-increasing, and strictly when it is strictly decreasing. For $k=0$, the non-strict conclusion still holds; the strict conclusion holds if $p(\cdot;t)\ne1+z^{n+1}$ for every $t\in J$.
\end{corollary}

\begin{proof}
Set $N=n+1$. The polynomial $p$ has constant and leading coefficients equal to $1$, and has simple zeros on $\T$. Its conjugate pair nearest to $1$ is $-\zeta,-\overline\zeta$. This is the unique pair of smallest absolute argument. By \cite[Theorem~1.1 and Remark~1]{BDPW91},
\[
 s(z)=\frac{p(z)}{(z+\zeta)(z+\overline\zeta)}
       =\sum_{l=0}^{N-2}s_lz^l
\]
has strictly positive coefficients. Write $p(z)=\sum_{j=0}^Np_jz^j$. Since $R(z)=P_{n+1}(z)/((z-\zeta)(z-\overline\zeta))=(-1)^Ns(-z)$, formula~\eqref{eq:A} gives
\begin{equation}\label{eq:extremesign}
 (-1)^{k+1}d_k=
 \sum_{\substack{0\leq j\leq N,\ 0\leq l\leq N-2\\[7pt]1+l-j=k}}p_js_l.
\end{equation}
For $1\leq k\leq n$, the term $j=0$, $l=k-1$ is positive. For $k=0$, the sum is positive precisely when at least one interior coefficient $p_1,\ldots,p_{N-1}$ is positive. Formula~\eqref{eq:conjugate} proves the stated conclusions.

It remains to verify the sufficient matrix condition. Expanding \eqref{eq:opuc} along its last row gives
\[
 Q_n(z)=\sum_{l=0}^na_lz^l,\quad
 (-1)^{n+l}a_l=\frac{\det C_{n+1}^{(n,l)}}{\det C_n}\geq0.
\]
The formula $P_{n+1}=zQ_n-(-1)^nQ_n^*$ now shows that every coefficient of $(-1)^{n+1}P_{n+1}(-z)$ is non-negative, as required.
\end{proof}

\begin{remark}\label{rem:finitepf}
For the sufficient condition in Corollary~\ref{cor:pf}, it is enough that the $n+1$ deleted minors $\det C_{n+1}^{(n,l)}$ be non-negative. No condition on minors involving moments of order greater than $n$ is needed. If only $c_0$ varies and the data are totally non-negative, the exceptional stationary case is $c_1=\cdots=c_n=0$: then $Q_n=z^n$ and $p=1+z^{n+1}$. In particular, strict total positivity is sufficient for the strict conclusion, but is not required.
\end{remark}

\subsection{Finite comparison of paraorthogonal zeros}

The preceding sign identities also permit direct comparison between two sets of finite moments. The semicircle hypothesis is required at only one endpoint; no assumption is imposed on the zeros along the interpolating segment.

\Needspace{9\baselineskip}
\begin{theorem}[One-endpoint comparison of zeros]\label{thm:zeroendpoints}
Let $n\geq1$, and let $A=(a_{r-s})_{r,s=0}^n$ and $B=(b_{r-s})_{r,s=0}^n$ be real symmetric positive-definite Toeplitz matrices. Let $Q_n^A$ and $Q_n^B$ be their monic orthogonal polynomials, defined by \eqref{eq:opuc}, and put
\[
 P_X(z;\epsilon)=zQ_n^X(z)-\epsilon(Q_n^X)^*(z),
 \quad X\in\{A,B\},\quad\epsilon\in\{-1,1\}.
\]
For the parameter specified in each case, write the arguments of the upper-half-plane zeros in increasing order as
\[
 0<\varphi_1^X<\cdots<\varphi_r^X<\pi,
 \quad r=\lfloor(n+1)/2\rfloor.
\]
\begin{enumerate}
\item[(i)] Put $\epsilon=(-1)^n$. If every zero of $P_B(\cdot;\epsilon)$ has positive real part and
\[
 (-1)^k(b_k-a_k)\geq0,\quad 0\leq k\leq n,
\]
then every zero of $P_A(\cdot;\epsilon)$ also has positive real part and
\[
 \varphi_j^A\leq\varphi_j^B,\quad 1\leq j\leq r.
\]
\item[(ii)] Put $\epsilon=-1$. If every zero of $P_B(\cdot;-1)$ has negative real part and
\[
 b_k-a_k\geq0,\quad 0\leq k\leq n,
\]
then every zero of $P_A(\cdot;-1)$ also has negative real part and
\[
 \varphi_j^A\geq\varphi_j^B,\quad 1\leq j\leq r.
\]
\end{enumerate}
In either case, all the angular inequalities are strict if $A\ne B$.
\end{theorem}

\begin{proof}
Consider $C(s)=(1-s)B+sA$, $0\leq s\leq1$. This segment is uniformly positive definite and extends to an open interval containing $[0,1]$ on which it remains positive definite. Theorem~\ref{thm:variation} therefore gives $C^1$ zero branches throughout the segment. The coefficients are real, and an upper zero cannot reach $1$ or $-1$: it would then coincide with its conjugate, contrary to simplicity. Distinct upper branches likewise cannot cross. When $n$ is even, the remaining real zero is fixed at $1$ in~\textup{(i)} and at $-1$ in~\textup{(ii)}; when $n$ is odd there are no real zeros for these parameters.

For~\textup{(i)}, write $\Delta_k=(-1)^k(b_k-a_k)\geq0$. As long as all zeros remain in the right semicircle, the coefficient calculation in the proof of Corollary~\ref{cor:right} gives $(-1)^{k+1}d_k=|d_k|>0$. Applying \eqref{eq:D} and \eqref{eq:conjugate} to an upper zero $\zeta=e^{i\varphi}$ yields
\[
 (h_\zeta+h_{\overline\zeta})\varphi'
 =-2\sin\varphi\sum_{k=0}^n(2-\delta_{k0})|d_k|\Delta_k\leq0.
\]
This inequality prevents any first exit from the right semicircle: no upper argument exceeds its initial value, which is strictly less than $\pi/2$. The real zeros and the lower zeros have already been accounted for. Thus the entire segment remains in the right semicircle. If $A\ne B$, at least one $\Delta_k$ is positive, so every upper angular derivative is strictly negative. Integration proves~\textup{(i)} and its strictness assertion.

For~\textup{(ii)}, put $\Delta_k=b_k-a_k\geq0$. The proof of Corollary~\ref{cor:left} gives $d_k>0$, and the corresponding identity is
\[
 (h_\zeta+h_{\overline\zeta})\varphi'
 =2\sin\varphi\sum_{k=0}^n(2-\delta_{k0})d_k\Delta_k\geq0.
\]
No upper argument falls below its initial value, which exceeds $\pi/2$. The same continuation and integration argument proves~\textup{(ii)}, again strictly when $A\ne B$.
\end{proof}

The following exact example illustrates simultaneous comparison. Take
\[
 A=\begin{pmatrix}13&5&-9\\[7pt]5&13&5\\[7pt]-9&5&13\end{pmatrix},
 \quad
 B=\begin{pmatrix}14&4&-7\\[7pt]4&14&4\\[7pt]-7&4&14\end{pmatrix}.
\]
Their leading principal determinants are $(13,144,44)$ and $(14,180,1386)$, respectively, and the signed moment differences are $(1,1,2)$. For $\epsilon=1$, direct substitution into \eqref{eq:opuc} gives
\[
 P_A(z;1)=(z-1)\left(z^2-\frac34z+1\right),\quad
 P_B(z;1)=(z-1)\left(z^2-\frac1{10}z+1\right).
\]
Thus $\det A=44<1386=\det B$ and $\varphi_1^A=\arccos(3/8)<\arccos(1/20)=\varphi_1^B$. The difference $B-A$ has the eigenvalue $-1$, with eigenvector $(1,0,-1)^\mathsf T$, so these comparisons do not follow from positive-semidefinite ordering of the matrices.

\subsection{A Rogers--Szeg\H{o} example and the limit of the conclusion}

\begin{example}\label{ex:rogers}
Fix $q\in(0,1)$ and $n\geq1$. Consider the rotated Rogers--Szeg\H{o} moments
\begin{equation}\label{eq:rogersmoments}
 c_j(0)=(-1)^j q^{j^2/2},\quad |j|\leq n.
\end{equation}
These data are positive definite. Indeed, diagonal sign changes followed by row and column factorisation reduce each leading principal determinant to a Vandermonde determinant:
\begin{equation}\label{eq:rogersdet}
\begin{aligned}
 \det C_M(0)
 &=q^{\sum_{r=0}^{M-1}r^2}
   \prod_{0\leq r<s\leq M-1}(q^{-s}-q^{-r})\\[7pt]
 &=\prod_{d=1}^{M-1}(1-q^d)^{M-d}>0,
 \quad 1\leq M\leq n+1.
\end{aligned}
\end{equation}
Indeed, after extracting $q^{r^2/2}$ and $q^{s^2/2}$ from row $r$ and column $s$, respectively, the remaining matrix is $(q^{-rs})$. Factoring $q^{-s}$ from each difference in its Vandermonde determinant cancels the prefactor. Empty products are understood to equal $1$. The monic orthogonal polynomial is
\begin{equation}\label{eq:rogersQ}
 Q_n(z;0)=\sum_{j=0}^n \qbinom{n}{j} q^{(n-j)/2}z^j,\quad
 \qbinom{n}{j}=\frac{(q;q)_n}{(q;q)_j(q;q)_{n-j}},
\end{equation}
where $(q;q)_m=\prod_{l=1}^m(1-q^l)$ and $(q;q)_0=1$; see \cite[Section~1.6]{S05I}. In particular, every coefficient is positive.

The normalisation in \eqref{eq:rogersQ} can also be checked directly. The finite $q$-binomial identity gives, for $0\leq\ell<n$,
\[
 \sum_{j=0}^n\qbinom{n}{j} q^{(n-j)/2}c_{j-\ell}(0)
 =(-1)^\ell q^{(n+\ell^2)/2}
   \prod_{r=0}^{n-1}(1-q^{r-\ell})=0,
\]
which is precisely the required orthogonality, since the moments are real and symmetric.

Choose $0\leq k\leq n$ and set $c_{\pm k}(t)=c_k(0)+u(t)$, where $u$ is a real $C^1$ function with $u(0)=0$, keeping the remaining finite moments fixed. Positive definiteness and positivity of the coefficients persist on a sufficiently small interval about $0$. With $b=-1$, all coefficients of $P_{n+1}(\cdot,-1;t)$ are therefore positive. If $\zeta(t)$ is its upper-half-plane zero nearest to $1$, \cite[Theorem~1.1]{BDPW91} implies positivity of every coefficient of
\[
 R(z;t)=\frac{P_{n+1}(z,-1;t)}{(z-\zeta(t))(z-\overline{\zeta(t)})}.
\]
Thus every Laurent coefficient of $z\overline{P_{n+1}(z,-1;t)}R(z;t)$ in the range $[-n,n]$ is positive. By \eqref{eq:conjugate}, $\zeta$ moves clockwise when $u$ increases, strictly when $u$ increases strictly. This conclusion holds for both parities and also permits a perturbation of $c_0$.

The other upper zeros need not move in the same direction. An exact example is obtained by taking
\begin{equation}\label{eq:counterdata}
 n=4,\quad q=\frac1{16},\quad k=2,\quad u(t)=t.
\end{equation}
Hence $c_{\pm2}(t)=1/256+t$, and all remaining moments in \eqref{eq:rogersmoments} remain fixed. To compute the perturbation explicitly, write
\begin{gather*}
 Q_4(z;t)=z^4+\alpha_3(t)z^3+\alpha_2(t)z^2+\alpha_1(t)z+\alpha_0(t),
 \\[7pt]
 v=\alpha_3+\alpha_0,\quad w=\alpha_2+\alpha_1.
\end{gather*}
The polynomial with $b=-1$ is then
\begin{equation}\label{eq:countervw}
 P_5(z;t)=z^5+v(t)z^4+w(t)(z^3+z^2)+v(t)z+1.
\end{equation}
Adding the first and fourth orthogonality equations, and then the second and third, gives the two-dimensional system
\begin{equation}\label{eq:countersystem}
 \begin{pmatrix}
  1+c_3&c_1+c_2(t)\\[7pt]
  c_1+c_2(t)&1+c_1
 \end{pmatrix}
 \begin{pmatrix}v(t)\\[7pt]w(t)\end{pmatrix}
 =-\begin{pmatrix}c_4+c_1\\[7pt]c_3+c_2(t)\end{pmatrix}.
\end{equation}
At $t=0$, put
\[
 a=\frac{262143}{262144},\quad d=\frac34,\quad
 e=-\frac{63}{256},\quad
 \Delta=ad-e^2=\frac{722925}{1048576}>0.
\]
Solving \eqref{eq:countersystem} gives
\[
 v(0)=\frac{4433}{16384},\quad
 w(0)=\frac{87637}{1048576}.
\]
Differentiating this system yields
\[
 av'(0)+ew'(0)=-w(0),\quad
 ev'(0)+dw'(0)=-1-v(0).
\]
Consequently,
\begin{align*}
 v'(0)&=\frac{-dw(0)+e(1+v(0))}{\Delta}
       =-\frac{257\cdot2091008}{987033600},\\[7pt]
 w'(0)&=\frac{ew(0)-a(1+v(0))}{\Delta}
       =-\frac{257\cdot7192437}{987033600}.
\end{align*}
Substitution into \eqref{eq:countervw} now gives
\begin{align}
 P_5(z;0)&=z^5+\frac{4433}{16384}z^4
 +\frac{87637}{1048576}(z^3+z^2)
 +\frac{4433}{16384}z+1,\label{eq:Pcounter}\\[7pt]
 \partial_tP_5(z;0)&=
 -\frac{257z(z+1)}{987033600}
 \bigl(2091008z^2+5101429z+2091008\bigr).\label{eq:Pcounterdot}
\end{align}
The reciprocal polynomial \eqref{eq:countervw} has the fixed zero $-1$. After dividing by $(z+1)z^2$ and putting $x=(z+z^{-1})/2$, the remaining zeros satisfy
\[
 F(x;t)=4x^2+2(v(t)-1)x+w(t)-v(t)-1=0.
\]
Equations~\eqref{eq:Pcounter}--\eqref{eq:Pcounterdot} yield
\begin{align}
 F(x;0)&=4x^2-\frac{11951}{8192}x-\frac{1244651}{1048576},\label{eq:Fcounter}\\[7pt]
 \partial_tF(x;0)&=-\frac{257(4182016x+5101429)}{987033600}<0,
 \quad -1\leq x\leq1.\label{eq:Fcounterdot}
\end{align}
The two roots satisfy $x_-\in(-2/5,-3/8)$ and $x_+\in(3/4,4/5)$. Let $\varphi_\pm(t)$ denote the arguments of the corresponding upper zero branches, with $\varphi_\pm(0)=\arccos x_\pm$. Subscripts on $F$ denote partial derivatives. Since $F_x(x_-;0)<0<F_x(x_+;0)$, implicit differentiation of $F(\cos\varphi;t)=0$ gives
\begin{equation}\label{eq:countervelocities}
 \varphi_-'(0)=\frac{F_t(x_-;0)}{\sqrt{1-x_-^2}F_x(x_-;0)}>0,
 \quad
 \varphi_+'(0)=\frac{F_t(x_+;0)}{\sqrt{1-x_+^2}F_x(x_+;0)}<0.
\end{equation}
Thus one upper zero moves anticlockwise and the other clockwise, although $c_2'=1$. Numerically,
\begin{align*}
 (\varphi_+(0),\varphi_+'(0))&\simeq(0.712368,-0.716541),\\[7pt]
 (\varphi_-(0),\varphi_-'(0))&\simeq(1.973710,0.213199).
\end{align*}
These decimal values illustrate the exact sign argument above. Figure~\ref{fig:counterexample} shows the zeros and their directions of motion. Thus the nearest-zero conclusion cannot be extended to all upper-half-plane zeros.
\end{example}

\begin{figure}[tb]
 \centering
 \setlength{\fboxsep}{4pt}
 \setlength{\fboxrule}{0pt}
 \resizebox{.76\textwidth}{!}{%
 \fbox{%
 \begin{tikzpicture}[x=2.2cm,y=2.2cm,font=\fontsize{9}{11}\selectfont,
   motion/.style={line width=.55pt,-{Stealth[length=4pt,width=3.2pt]}}]
   \pgfmathsetmacro{\angleplus}{acos(0.7568160422462677)}
   \pgfmathsetmacro{\angleminus}{acos(-0.3921004660743927)}
   \draw[gray!45,line width=.3pt] (-1.12,0)--(1.12,0);
   \draw[gray!45,line width=.3pt] (0,-1.12)--(0,1.12);
   \draw[gray,line width=.45pt] (0,0) circle[radius=1];
   \fill (-1,0) circle[radius=1.35pt];
   \node[anchor=east] at (-1.08,-.06) {$-1$ (fixed)};
   \node[anchor=west] at (1.08,-.04) {$1$};
   \node[anchor=west] at (.025,1.10) {$i$};
   \node[anchor=west] at (.025,-1.10) {$-i$};
   \foreach \a/\b in {\angleplus/{\angleplus-11},
     -\angleplus/{-\angleplus+11},
     \angleminus/{\angleminus+11},
     -\angleminus/{-\angleminus-11}} {
     \fill (\a:1) circle[radius=1.35pt];
     \draw[motion] (\a:1.025) arc[start angle=\a,end angle=\b,radius=1.025];
   }
   \node[anchor=west] at (.98,.88) {$\varphi_+'(0)<0$};
   \node[anchor=east] at (-.70,1.09) {$\varphi_-'(0)>0$};
   \useasboundingbox (-1.9,-1.26) rectangle (1.88,1.32);
 \end{tikzpicture}%
 }%
 }
 \caption{Zeros of $P_5(\cdot,-1;0)$ for the data in \eqref{eq:counterdata}. Arrows show the directions when $c_2$ increases; the zero $-1$ is fixed. The two upper zeros move in opposite directions.}
 \label{fig:counterexample}
\end{figure}

\section{Proofs of the main results}\label{sec:proofs}

\subsection{Finite comparison and strictness}

\begin{proof}[Proof of Theorem~\ref{thm:comparison}]
Assume first~\textup{(i)}, and let $A$ satisfy the hypotheses in~\textup{(ii)}. Put $N=n+1$ and $\Delta_k=(-1)^k(b_k-a_k)$. Grouping the inverse/cofactor formula by diagonals gives
\begin{equation}\label{eq:traceorder}
 \tau:=\tr\bigl(B^{-1}(B-A)\bigr)
 =\frac{1}{\det B}\sum_{k=0}^n\Delta_kS_k(B)\geq0.
\end{equation}
Here $(-1)^{r+s}=(-1)^{|r-s|}$, so all cofactors on the $k$th pair of diagonals have the common sign factor $(-1)^k$. Put $X=B^{-1/2}AB^{-1/2}>0$. Then \eqref{eq:traceorder} gives $\tr X=N-\tau\leq N$. Applying the arithmetic--geometric mean inequality to the positive eigenvalues of $X$ yields
\[
 \frac{\det A}{\det B}=\det X
 \leq\left(\frac{\tr X}{N}\right)^N\leq1.
\]
Equality forces all eigenvalues of $X$ to equal $1$. Because $X$ is Hermitian, this means $X=I$, or $A=B$. The converse is immediate.

To prove necessity, assume~\textup{(ii)} and fix $0\leq k\leq n$. Let $H^{[k]}$ be the real symmetric Toeplitz matrix of order $n+1$ with entries
\[
 H^{[k]}_{rs}=\begin{cases}
 (-1)^k,& |r-s|=k,\\[7pt]
 0,& |r-s|\ne k.
 \end{cases}
\]
For all sufficiently small $u>0$, the matrix $A_u=B-uH^{[k]}$ is positive definite and satisfies \eqref{eq:momentorder}. Thus $\log\det A_u\leq\log\det B$, with equality at $u=0$, and differentiation at $0$ from the right gives
\[
 0\geq\left.\frac{\dd}{\dd u}\log\det A_u\right|_{u=0}
 =-\tr(B^{-1}H^{[k]})=-\frac{S_k(B)}{\det B}.
\]
This proves~\textup{(i)}. Finally, $S_0(B)>0$ because every principal deleted minor of a positive-definite matrix is positive.
\end{proof}

The argument also yields a quantitative estimate. With $X$, $N$ and $\tau$ as above, and with $\|\cdot\|_{\mathrm F}$ and $\|\cdot\|_2$ denoting the Frobenius and spectral norms, respectively,
\begin{equation}\label{eq:detgap}
 \log\frac{\det B}{\det A}
 \geq\tau+\frac{\|X-I\|_{\mathrm F}^2}{2\max\{1,\|X\|_2\}}
 \geq\tau+\frac{\|X-I\|_{\mathrm F}^2}{2N}.
\end{equation}
Indeed, integrating $(t-1)/t$ between $1$ and $x>0$ gives
\[
 x-1-\log x\geq\frac{(x-1)^2}{2\max\{1,x\}}.
\]
Summing over the eigenvalues of $X$ and using $N-\tr X=\tau$ and $\|X\|_2\leq\tr X\leq N$ proves the estimate. A small logarithmic determinant gap therefore controls the distance of the normalised matrix $X$ from the identity, even when $\tau=0$.

For comparison with the classical differential argument, Jacobi's formula gives, for real $C^1$ moments,
\begin{equation}\label{eq:jacobi}
 \frac{\dd}{\dd t}\det C_{n+1}
 =\tr(\adj(C_{n+1})C_{n+1}')
 =\sum_{r,s=0}^n(-1)^{r+s}\det C_{n+1}^{(r,s)}c_{r-s}'.
\end{equation}
If only $c_{\pm k}$ vary, this reduces to
\[
 \frac{\dd}{\dd t}\det C_{n+1}
 =(-1)^kc_k'\sum_{|r-s|=k}\det C_{n+1}^{(r,s)}.
\]
The sum includes both diagonals when $k>0$ and the main diagonal once when $k=0$. This accounts for the non-strict derivative sign. The endpoint argument establishes strict comparison without assuming positivity of every deleted minor, even when the derivative of a strictly monotone moment vanishes at some points.

\subsection{The half-plane criterion}

We begin with the CD identity at a paraorthogonal zero. With the conventions \eqref{eq:kernel}--\eqref{eq:popuc}, any zero $\xi$ of $P=P_{n+1}(\cdot,b)$ satisfies
\begin{equation}\label{eq:CDzero}
 K_n(\xi,z)=\kappa_n^2\overline{Q_n(\xi)}\,
 \frac{P(z)}{z-\xi}.
\end{equation}
Indeed, the usual CD formula reads
\[
 K_n(w,z)=\kappa_n^2
 \frac{\overline{Q_n^*(w)}Q_n^*(z)
       -\overline w z\,\overline{Q_n(w)}Q_n(z)}{1-\overline w z}.
\]
For $w=\xi$, substitute $Q_n^*(\xi)=b\xi Q_n(\xi)$ and cancel the factor $-\overline\xi$. The identity extends to $z=\xi$ by continuity. In particular, the normalised cardinal polynomials
\begin{equation}\label{eq:cardinal}
 \ell_m(z)=\frac{P(z)}{(z-\zeta_m)P'(\zeta_m)}
 =\frac{K_n(\zeta_m,z)}{K_n(\zeta_m,\zeta_m)}
\end{equation}
at the $N=n+1$ distinct zeros $\zeta_1,\ldots,\zeta_N$ are mutually orthogonal; here $P'$ denotes differentiation with respect to $z$. They satisfy
\begin{equation}\label{eq:weights}
 \langle\ell_m,\ell_l\rangle=\lambda_m\delta_{ml},\quad
 \lambda_m=\frac1{K_n(\zeta_m,\zeta_m)}>0.
\end{equation}
These identities also give the finite positive quadrature construction, although no integral representation is needed here.

\begin{proof}[Proof of Theorem~\ref{thm:hurwitz}]
Suppress the moment data from the notation, and write $C=C_{n+1}$ and $N=n+1$. The real polynomial $P=P_{n+1}(\cdot,(-1)^n)$ has distinct zeros $\zeta_1,\ldots,\zeta_N$ in $\T_+$. Define
\[
 R_m(z)=\frac{P(z)}{z-\zeta_m}=\sum_{j=0}^nd_{mj}z^j.
\]
Let $V$ have as its $m$th column the coefficients of $\ell_m=R_m/P'(\zeta_m)$. These cardinal polynomials form a basis of $\Pcal_n$. Thus \eqref{eq:weights} says $V^*CV=\operatorname{diag}(\lambda_m)$, and inversion gives
\begin{equation}\label{eq:inversegram}
 (C^{-1})_{rs}=\sum_{m=1}^N
 \frac{d_{mr}\overline{d_{ms}}}{\lambda_m|P'(\zeta_m)|^2}.
\end{equation}
Put $\varepsilon_j=(-1)^{n-j}$. We claim that, for every $m,r,s$,
\begin{equation}\label{eq:coeffpositive}
 \Re\bigl((\varepsilon_rd_{mr})\overline{\varepsilon_sd_{ms}}\bigr)>0.
\end{equation}
If $\zeta_m=1$, $R_m$ is real monic with roots in the open right half-plane. Therefore $(-1)^nR_m(-z)$ has strictly positive coefficients, and $\varepsilon_jd_{mj}>0$. This proves \eqref{eq:coeffpositive} in this case.

Otherwise $\zeta_m$ is non-real, and
\[
 R_m(z)=(z-\overline{\zeta_m})H_m(z),\quad
 H_m(z)=\sum_{j=0}^{n-1}f_jz^j,
\]
where $H_m$ is real monic and again has all its zeros in the open right half-plane. Set $f_{-1}=f_n=0$. The polynomial $(-1)^{n-1}H_m(-z)$ has strictly positive coefficients, so
\[
 \varepsilon_jd_{mj}=A_j+B_j\overline{\zeta_m},\quad
 A_j=\varepsilon_jf_{j-1}\geq0,\quad
 B_j=-\varepsilon_jf_j\geq0,
\]
and $(A_j,B_j)\ne(0,0)$ for every $0\leq j\leq n$. Consequently the real part in \eqref{eq:coeffpositive} equals
\[
 A_rA_s+B_rB_s+(A_rB_s+B_rA_s)\Re\zeta_m>0.
\]
At least one of the four products is positive, and $\Re\zeta_m>0$.

Since $C^{-1}$ is real, \eqref{eq:inversegram} now implies
\[
 (-1)^{r+s}(C^{-1})_{rs}>0.
\]
Finally, $\det C^{(r,s)}=(-1)^{r+s}\det C\,(C^{-1})_{sr}>0$. This proves the theorem in both parities.
\end{proof}

\subsection{Variation through the moment functional}

\begin{proof}[Proof of Theorem~\ref{thm:variation}]
The coefficients of $Q_n$ are rational functions of the finite moments, with the positive determinant $\det C_n$ in the denominator. Hence the coefficients of $P$ are $C^1$. Its zeros are simple and remain on $\T$, so the implicit function theorem gives local $C^1$ zero branches. They continue uniquely across every compact subinterval of $J$: the zeros cannot collide, and the unit circle is compact. This gives global branches on $J$, each with a real $C^1$ argument because $J$ is an interval.

For a zero $\xi$, put $R_\xi(z)=P(z)/(z-\xi)$. The reproducing identity and \eqref{eq:CDzero} give, for every $f\in\Pcal_n$,
\begin{equation}\label{eq:reproduction}
 \Lcal_t[\overline{R_\xi}f]
 =\frac{f(\xi)}{\kappa_n^2Q_n(\xi)},\quad
 h_\xi=\frac{P'(\xi)}{\kappa_n^2Q_n(\xi)}>0.
\end{equation}
Positivity in the second identity also follows directly from \eqref{eq:h}. Since $P$ is monic, $\partial_tP\in\Pcal_n$. Differentiating $P(\xi(t);t)=0$ and using \eqref{eq:reproduction}, with $\xi=e^{i\theta_\xi}$, yields
\begin{equation}\label{eq:onezero}
 h_\xi\theta_\xi'
 =i\overline\xi\,\Lcal_t[\overline{R_\xi}\,\partial_tP]
 =-i\Lcal_t\!\left[\frac{z\overline P\,\partial_tP}{z-\xi}\right].
\end{equation}
The last equality uses $\overline{z-\xi}=-(z-\xi)/(z\xi)$ on $\T$.

We first prove the individual identity. Differentiate $\Lcal_t[\overline{q_j}Q_n]=0$ for $0\leq j<n$. The contribution from $\partial_tq_j$ vanishes by orthogonality. Since $\partial_tQ_n\in\Pcal_{n-1}$, expansion in the orthonormal basis gives
\begin{equation}\label{eq:Qvariation}
 \partial_tQ_n(\xi;t)
 =-\dot\Lcal_t[Q_n(z;t)K_{n-1}(z,\xi;t)].
\end{equation}
At a zero $\xi=e^{i\theta_\xi}$, the defining equation is
\[
 \xi\frac{Q_n(\xi;t)}{Q_n^*(\xi;t)}=e^{-i\beta(t)}.
\]
On $\T$, $Q_n^*(e^{i\theta};t)=e^{in\theta}\overline{Q_n(e^{i\theta};t)}$. Differentiating the argument of the preceding equation along the zero branch therefore gives
\[
 \left(1-n+2\Re\frac{\xi Q_n'(\xi;t)}{Q_n(\xi;t)}\right)\theta_\xi'
 =-\beta'-2\Im\frac{\partial_tQ_n(\xi;t)}{Q_n(\xi;t)}.
\]
The factor in parentheses equals $P'(\xi;t)/Q_n(\xi;t)=h_\xi/\nu$, by differentiation of $P=zQ_n-\overline bQ_n^*$ and \eqref{eq:reproduction}. Multiplication by $\nu$ and use of \eqref{eq:Qvariation} prove \eqref{eq:individual}, since a Hermitian functional commutes with taking real and imaginary parts. The expression in \eqref{eq:V} belongs to $\Lambda_n$: before taking its imaginary part, its Laurent exponents range from $-(n-1)$ to $n$.

The orthogonality relations for $Q_n$ and $Q_n^*$ imply
\begin{equation}\label{eq:quasi}
 \Lcal_t[P\overline g]=0,\quad g\in\Pcal_n,\quad g(0)=0.
\end{equation}
Indeed, it suffices to take $g=z^j$, $1\leq j\leq n$: the term $zQ_n$ vanishes by orthogonality against $z^{j-1}$, and the term $Q_n^*$ by conjugating orthogonality against $z^{n-j}$. Every product in \eqref{eq:quasi} belongs to $\Lambda_n$.

For the two distinct branches $\zeta,\eta$, let
\[
 H(z;t)=\frac{zP(z;t)}{(z-\zeta(t))(z-\eta(t))}.
\]
This is a $C^1$ polynomial of degree $n$ with $H(0;t)=0$. Subtracting \eqref{eq:onezero} at the two zeros, and using
\[
 \overline H=\zeta\eta\,
 \frac{z\overline P}{(z-\zeta)(z-\eta)}\quad\hbox{on }\T,
\]
gives
\begin{equation}\label{eq:subtract}
 h_\zeta\varphi'-h_\eta\psi'
 =-\frac{i(\zeta-\eta)}{\zeta\eta}\,
 \Lcal_t[\overline H\,\partial_tP].
\end{equation}
Differentiate \eqref{eq:quasi} with $g=H$. Since $\partial_tH$ has degree at most $n$ and vanishes at $0$, its contribution vanishes by the same relation. Hence
\[
 \Lcal_t[\overline H\,\partial_tP]=-\dot\Lcal_t[P\overline H].
\]
Substitution into \eqref{eq:subtract} proves \eqref{eq:twozero}.

The Laurent polynomial $B_{\zeta,\eta}=i(\zeta-\eta)H\overline P$ lies in $\Lambda_n$ because $H(0)=0$ and $\deg H=n$. It is real on $\T$: conjugation of $z/((z-\zeta)(z-\eta))$ multiplies it by $\zeta\eta$, whereas conjugation of $i(\zeta-\eta)$ divides it by $\zeta\eta$. Thus its coefficients satisfy $f_{-j}=\overline{f_j}$. Applying $\dot\Lcal_t[z^j]=c_{-j}'$ and pairing $j$ with $-j$ gives \eqref{eq:D}. Every differentiation above takes place in a fixed finite-dimensional vector space. No representing measure is differentiated.
\end{proof}

It remains to verify \eqref{eq:bkconj} and \eqref{eq:bkfixed}. For any two distinct zeros $\zeta,\eta$, \eqref{eq:CDzero} and its conjugate give, on $\T$,
\begin{equation}\label{eq:kernelproduct}
 K_n(\zeta,z)K_n(z,\eta)
 =-\kappa_n^4\eta\,\overline{Q_n(\zeta)}Q_n(\eta)
 \frac{z|P(z)|^2}{(z-\zeta)(z-\eta)}.
\end{equation}
For $\eta=\overline\zeta$, multiplication by $-\zeta Q_n(\zeta)/Q_n(\overline\zeta)$ gives $\kappa_n^4|Q_n(\zeta)|^2A_\zeta$. Extracting the coefficient of $z^k$ proves the second equality in \eqref{eq:bkconj}. The first follows from
\[
 A_\zeta(e^{i\theta})
 =-\frac{|K_n(\zeta,e^{i\theta})|^2}{\kappa_n^4|Q_n(\zeta)|^2}
 \frac{\sin((\theta-\varphi)/2)}{\sin((\theta+\varphi)/2)}.
\]
For a fixed $\eta$, multiplication of \eqref{eq:kernelproduct} by
$-i\overline\eta(\zeta-\eta)/(Q_n(\eta)\overline{Q_n(\zeta)})$
gives $\kappa_n^4B_{\zeta,\eta}$ and proves \eqref{eq:bkfixed}. These calculations also show directly that all apparent singularities cancel.

\section*{Acknowledgements}
\begingroup
\emergencystretch=1em

The first author acknowledges financial support from the Centre for Mathematics of the University of Coimbra (CMUC), funded by the Portuguese Foundation for Science and Technology (FCT), under the projects UID/00324/2025 (\url{https://doi.org/10.54499/UID/00324/2025}) and UID/PRR/00324/2025. The first author also acknowledges financial support from the FCT under the grant \url{https://doi.org/10.54499/2022.00143.CEECIND/CP1714/CT0002}. The second author was supported by FCT grant 2021.05089.BD.

\par\endgroup
\makeatletter
\renewcommand{\@biblabel}[1]{#1.}
\makeatother
\providecommand{\bysame}{\leavevmode\hbox to3em{\hrulefill}\thinspace}

\end{document}